\documentclass[12pt,reqno]{amsart}
\usepackage{amssymb}
\usepackage{}
\usepackage{amsfonts}
\usepackage{mathrsfs}
\usepackage{amsmath,amsthm,amssymb,amsfonts,amscd}
\usepackage{mathrsfs}
\usepackage{bbding}
\usepackage{graphicx,latexsym}

\theoremstyle{plain}
\newtheorem{theorem}{Theorem}[section]

\newtheorem{lemma}{Lemma}[section]

\theoremstyle{definition}

\newtheorem{remark}[lemma]{Remark}

\def \N {{\mathbb N}}
\def \Z {{\mathbb Z}}

\def\bg{\bigg}
\def\({\bg(}
\def\){\bg)}

\def\le {\leqslant}
\def\ge {\geqslant}

\begin{document}
\medskip

\title{On a discriminator for the polynomial $f(x)=x^3+x$}

\author{Quan-Hui Yang}
\address{School of Mathematics and Statistics \\ Nanjing University of Information Science and Technology \\ Nanjing 210044 \\Reople's Republic of China}
\email{yangquanhui01@163.com}
\author{Lilu Zhao}
\address{School of Mathematics \\ Shandong University \\ Jinan  250100 \\Reople's Republic of China}
\email{zhaolilu@sdu.edu.cn}

\begin{abstract} Let $\Delta(n)$ denote the smallest positive integer $m$ such that $a^3+a(1\le a\le n)$ are pairwise distinct modulo $m$. The purpose of this paper is to determine $\Delta(n)$ for all positive integers $n$.
\end{abstract}
\thanks{
\newline\indent This work is supported by National Natural Science Foundation of China (Grant No. 11922113).
}

\maketitle

\section{Introduction}
\setcounter{lemma}{0}
\setcounter{theorem}{0}
\setcounter{corollary}{0}
\setcounter{equation}{0}

\medskip

 For a polynomial $f(x)\in \Z[x]$ with all $f(a)(a\in \Z^{+})$ pairwise distinct, we introduce the discriminator $\Delta_f(n)$ defined to be the smallest positive integer $m$ such that $f(a)(1\le a\le n)$ are pairwise distinct modulo $m$.

As a simple application of Bertrand's postulate, Arnold, Benkoski and McCabe \cite{ABM} determined $\Delta_{f}(n)$ for $f(x)=x^2$, and they showed that for $n>4$, $\Delta_{f}(n)$ is the smallest positive integer $m\ge 2n$ such that $m$ is $p$ or $2p$ with $p$ an odd prime. Sun \cite{S13} studied $\Delta_f(n)$ for other quadratic polynomials. For example, it was proved in \cite{S13} that if $f(x)=2x(x-1)$ then $\Delta_{f}(n)$ is the least prime number greater than $2n-2$, and in particular $\Delta_{f}(n)$ runs over all prime values.

Among other things, Schumer \cite{S} studied $\Delta_f(n)$ with $f(x)=x^3$. For the study of discriminator $\Delta_f(n)$ with other higher degree polynomials $f$, one may refer to \cite{BSW,Moree,MM,Zieve}. In this paper, we focus on $\Delta_f(n)$ with $f(x)=x^3+x$.

The main result in this paper is the following.
\begin{theorem}\label{theorem1}Let $\Delta(n)=\Delta_f(n)$ with $f(x)=x^3+x$. We have
\begin{align*}\Delta(n)=\begin{cases}7\cdot 3^{6s+4} & \ \textrm{ if } n=3^{6s+5}+1 \textrm{ or } n=3^{6s+5}+2 \textrm{ for some } s\in \N,
\\
3^{\lceil \log_3 n\rceil} & \ \textrm{ otherwise},\end{cases}\end{align*}
where $\lceil x \rceil$ denotes the smallest integer no less than $x$.
\end{theorem}

A closely related problem is to determine $D(n)$, which denotes the smallest positive integer $m$ such that $a^3+a(1\le a\le n)$ are pairwise distinct modulo $m^2$. The authors \cite{YZ} proved that $D(n)=3^{\lceil \log_3 \sqrt{n}\rceil}$, which was conjectured by Z.-W. Sun (see Conjecture 6.76 in \cite{S2021}). The present work is motivated by the above original conjecture of Sun. Different from $D(n)$, the discriminator $\Delta(n)$ is not always a power of three. For example, $\Delta(n)=7\cdot 3^{4}$ when $n=244$ or $245$. This was first observed by Sun (see the remark to Conjecture 6.76 in \cite{S2021}). According to Theorem \ref{theorem1}, the third example of $n$ satisfying $\Delta(n)\not=3^{\lceil \log_3 n\rceil}$ is over $10^5$.

We prove Theorem \ref{theorem1} by combining methods from elementary number theory and analytic number theory. We point out that in order to deal with $\Delta(n)$ we have to study an incomplete character sum, which is not involved in the work \cite{YZ}. The incomplete character sum will be handed by the elementary method when the length of the summation is about $\frac{p}{4}$, and it will be handed by the analytic method when the length of the summation is about $\frac{p}{6}$. The details will be given in Section 3. Moreover, in the very special case $p=7$, we have to discuss the value of Legendre symbol separately (see Lemma \ref{lemma51} in Section 5).

\bigskip

We use the following notations in this paper. Let $\Z^{+}$ denote the set of all positive integers and let $\N=\Z^{+}\cup\{0\}$. We use $e(\alpha)$ to denote $e^{2\pi i\alpha}$. The notation $\lceil x \rceil$ denotes the smallest integer no less than $x$, and $\lfloor x\rfloor$ denotes the greatest integer no more than $x$.

\bigskip
\section{Preparations}

\setcounter{lemma}{0}
\setcounter{theorem}{0}
\setcounter{corollary}{0}
\setcounter{equation}{0}

We introduce
\begin{align*}\mathcal{E}=\{3^{6s+5}+1:\ s\in \N\} \cup \{ 3^{6s+5}+2:\ s\in \N\}.\end{align*}
Throughout this paper, we use the letter $k$ to denote
\begin{align*}k=\lceil \log_3 n\rceil.\end{align*}
We first point out that $a^3+a(1\le a\le n)$ are pairwise distinct modulo $3^{k}$, and therefore $n\le \Delta(n)\le 3^k$. In order to establish Theorem \ref{theorem1}, it suffices to prove the following two results.
\begin{lemma}\label{lemma21}Let $n\not\in \mathcal{E}$. Suppose that
\begin{align}\label{inequality1} n\le m<3^k<3n.\end{align}
Then there exist  $1\le a<b\le n$ such that $b^3+b\equiv a^3+a\pmod{m}$.\end{lemma}

\begin{lemma}\label{lemma22}Let $n=3^{6s+5}+1$ or $3^{6s+5}+2$ with $s\in \N$. Suppose that
\begin{align}\label{inequality2} n\le m<7\cdot 3^{6s+4}.\end{align}
Then there exist  $1\le a<b\le n$ such that $b^3+b\equiv a^3+a\pmod{m}$. Moreover, $a^3+a(1\le a\le n)$ are pairwise distinct modulo $7\cdot 3^{6s+4}$.
\end{lemma}

We shall consider the following $8$ cases.

(i) $m=\delta p$, where $\delta\ge 6$, $p\ge 5$ is a prime, $p\not=7$ and $p\nmid \delta$.

(ii) $m=\delta p^r$, where $\delta\ge 4$, $p\ge 5$ is a prime, $r\ge 2$ is a positive integer.

(iii) $m=2^r$, where $r\in \Z^{+}$.

(iv) $m=2^rt$, where $t\ge 5$ is an odd number and $r\ge 2$.

(v) $m=2^r3^s$, where $r,s\in \Z^{+}$.

(vi) $m=3^{r}\cdot 14$, where $r\in \N$.

(vii) $m=\delta p^{t}$, where $1\le \delta \le 3$, $p\ge 5$ is a prime, $t\in \Z^{+}$.

(viii) $m=3^{r}\cdot 7$, where $r\in \N$.

\noindent The letter $\delta$ always denotes a positive integer. Note that \eqref{inequality2} implies \eqref{inequality1}. Throughout this paper, we assume that \eqref{inequality1} holds.

\bigskip

\section{An incomplete character sum}

\setcounter{equation}{0}

\medskip

For $u\in \Z$, $\delta\in \Z^{+}$ and $p\ge 3$, we introduce
\begin{align*}A_p(\delta,u)=\sum_{-\frac{p-1}{2}\le x\le \frac{p-1}{2}}\Big(\frac{\delta^2x^2+4}{p}\Big)e\big(\frac{ux}{p}\big),\end{align*}
where $\big(\frac{\cdot}{p}\big)$ denotes the Legendre symbol.

\begin{lemma}\label{lemma31}Suppose that $p\ge 3$ is a prime and $p\nmid \delta$.

(i) If $p|u$, then $A_p(\delta,u)=-1$.

(ii) If $p\nmid u$, then $|A_p(\delta,u)|\le 2\sqrt{p}$.\end{lemma}
\begin{proof}For an odd prime $p$, it is well-known that
$$\sum_{1\le c\le p-1}\Big(\frac{c}{p}\Big)e(\frac{c}{p})=\sum_{1\le x\le p}e\big(\frac{x^2}{p}\big),$$
and $|\tau_p|=\sqrt{p}$, where $\tau_p$ denotes the above Gauss sum. By
$$\sum_{1\le c\le p-1}\Big(\frac{c}{p}\Big)e(\frac{c(\delta^2x^2+4)}{p})=\Big(\frac{\delta^2x^2+4}{p}\Big)\tau_p,$$
we deduce that
\begin{align*}A_p(\delta,u)=&\frac{1}{\tau_p}\sum_{-\frac{p-1}{2}\le x\le \frac{p-1}{2}}e(\frac{ux}{p})\sum_{1\le c\le p-1}\Big(\frac{c}{p}\Big)e(\frac{c(\delta^2x^2+4)}{p})
\\= & \frac{1}{\tau_p}\sum_{1\le c\le p-1}e(\frac{4c}{p})\Big(\frac{c}{p}\Big)\sum_{-\frac{p-1}{2}\le x\le \frac{p-1}{2}}e(\frac{c\delta^2x^2+ux)}{p}.\end{align*}
Note that
\begin{align*}\sum_{-\frac{p-1}{2}\le x\le \frac{p-1}{2}}e(\frac{c\delta^2x^2+ux)}{p}=\Big(\frac{c}{p}\Big)e(\frac{-\overline{4\delta^2c}\, u^2}{p})\tau_p,\end{align*}
where $\overline{d}$ means $\overline{d}\cdot d\equiv 1\pmod{p}$. Now we conclude that
\begin{align}\label{resultAp}A_p(\delta,u)=\sum_{1\le c\le p-1}e(\frac{-\overline{4\delta^2c}\, u^2+4c}{p}).\end{align}
If $p|u$, then the summation in \eqref{resultAp} is a Ramanujan sum and $A_p(\delta,u)=1$. If $p\nmid u$, then by Weil's bound on Kloosterman sums (see (4.19) in \cite{I}) we have $|A_p(\delta,u)|\le 2\sqrt{p}$. This completes the proof.
\end{proof}

We remark that Lemma \ref{lemma31} (i) is a well-known result. For a prime $p\ge 5$ and $p\nmid \delta$, we define $\ell_p(\delta)$ to be smallest positive integer $x$ such that
$$\Big(\frac{-3\delta^2x^2-12}{p}\Big)\in\{0,1\}.$$
We introduce
\begin{align*}L_p=\begin{cases}\frac{p+3}{4}\ \ \textrm{ if } p\equiv 1\pmod{12},
\\ \frac{p-1}{4}\ \ \textrm{ if } p\equiv 5\pmod{12},
\\ \frac{p+5}{4}\ \ \textrm{ if } p\equiv 7\pmod{12},
\\ \frac{p+1}{4}\ \ \textrm{ if } p\equiv 11\pmod{12}.\end{cases}\end{align*}

We point out that $L_p<\frac{p}{3}$ holds for $p\ge 5$ except $p=7$.
\begin{lemma}\label{lemma32}Suppose that $p\ge 5$ is a prime and $p\nmid \delta$. We have
$$\ell_p(\delta)\le L_p.$$\end{lemma}
\begin{proof}By Lemma \ref{lemma31} (i), we have
\begin{align*}A_p(\delta,0)=2\sum_{1\le x\le \frac{p-1}{2}}\Big(\frac{\delta^2x^2+4}{p}\Big)+1=-1,\end{align*}
and therefore,
\begin{align}\label{halfsum}\sum_{1\le x\le \frac{p-1}{2}}\Big(\frac{\delta^2x^2+4}{p}\Big)=-1.\end{align}
We introduce 
\begin{align*}N_p^{+}=&|\{1\le x\le \frac{p-1}{2}:\ \big(\frac{\delta^2x^2+4}{p}\big)=+1\}|,
\\ N_p^{-}=&|\{1\le x\le \frac{p-1}{2}:\ \big(\frac{\delta^2x^2+4}{p}\big)=-1\}|,
\\ N_p^{0}=&|\{1\le x\le \frac{p-1}{2}:\ \big(\frac{\delta^2x^2+4}{p}\big)=0\}|.\end{align*}

In view of \eqref{halfsum}, we have the following conclusions. If $p\equiv 1\pmod{4}$, then $N_p^{0}=1$, $N_p^{+}=\frac{p-5}{4}$ and $N_p^{-}=\frac{p-1}{4}$. If $p\equiv 3\pmod{4}$, then $N_p^{0}=0$, $N_p^{+}=\frac{p-3}{4}$ and $N_p^{-}=\frac{p+1}{4}$.

Case $p\equiv 1\pmod{12}$. We have $\big(\frac{-3}{p}\big)=1$ and $\ell_p(\delta)$ is the smallest positive integer $x$ such that
$\big(\frac{\delta^2x^2+4}{p}\big)\in\{0,1\}$. Note that $N_p^{0}+N_p^{+}=\frac{p-1}{4}$. Now we conclude that $\ell_p(\delta)\le \frac{p-1}{2}-(N_p^{0}+N_p^{+})+1=\frac{p-1}{2}-\frac{p-1}{4}+1=L_p$.

Case $p\equiv 5\pmod{12}$. We have $\big(\frac{-3}{p}\big)=-1$ and $\ell_p(\delta)$ is the smallest positive integer $x$ such that
$\big(\frac{\delta^2x^2+4}{p}\big)\in\{0,-1\}$. Note that $N_p^{0}+N_p^{-}=\frac{p+3}{4}$. Now we conclude that $\ell_p(\delta)\le \frac{p-1}{2}-(N_p^{0}+N_p^{-})+1=\frac{p-1}{2}-\frac{p+3}{4}+1=L_p$.

Case $p\equiv 7\pmod{12}$. We have $\big(\frac{-3}{p}\big)=1$ and $\ell_p(\delta)$ is the smallest positive integer $x$ such that
$\big(\frac{\delta^2x^2+4}{p}\big)=1$. Note that $N_p^{+}=\frac{p-3}{4}$. Now we conclude that $\ell_p(\delta)\le \frac{p-1}{2}-N_p^{+}+1=\frac{p-1}{2}-\frac{p-3}{4}+1=L_p$.

Case $p\equiv 11\pmod{12}$. We have $\big(\frac{-3}{p}\big)=-1$ and $\ell_p(\delta)$ is the smallest positive integer $x$ such that
$\big(\frac{\delta^2x^2+4}{p}\big)=-1$. Note that $N_p^{-}=\frac{p+1}{4}$. Now we conclude that $\ell_p(\delta)\le \frac{p-1}{2}-N_p^{-}+1= \frac{p-1}{2}-\frac{p+1}{4}+1=L_p$.

We are done.
\end{proof}

\begin{lemma}\label{lemma33}Suppose that $m=\delta p^{r}$, where $\delta, r\in \Z^{+}$, $p\ge 5$ is a prime and $p\nmid \delta$.

(i) If $p^{r}+\delta\frac{p-1}{2}\le n$, then there exist $1\le a<b\le n$ such that $b^3+b\equiv a^3+a\pmod{m}$.

(ii) If $r=1$ and $p+\delta \ell_p(\delta)\le n$, then there exist  $1\le a<b\le n$ such that $b^3+b\equiv a^3+a\pmod{m}$.

\end{lemma}
\begin{proof}We consider $1\le a\le p^{r}$ and $b=a+\delta c$ with $c\in \Z^{+}$. It suffices to find $a,c\in \Z^{+}$ such that
$a+\delta c\le n$ and
$$a^2+a(a+\delta c)+(a+\delta c)^2+1\equiv 0\pmod{p^{r}},$$
which is equivalent to
\begin{align}\label{congruencequad}(6a+3\delta c)^2\equiv -3\delta^2c^2-12\pmod{p^{r}}.\end{align}

In view of \eqref{halfsum}, we conclude that there exists $1\le c\le \frac{p-1}{2}$ such that $-3\delta^4c^2-12$ is a quadratic residue modulo $p$. Then it is easy to deduce that there exists $1\le a\le p^{r}$ such that $(6a+3\delta c)^2\equiv -3\delta^2c^2-12\pmod{p^{r}}$. This completes the proof the conclusion (i).

By the definition of $\ell_p(\delta)$, we can find $1\le c\le \ell_p(\delta)$ such that $\big(\frac{-3\delta^4c^2-12}{p}\big)\in \{0,1\}$. Then we can find $1\le a\le p$ such that $(6a+3\delta c)^2\equiv -3\delta^2c^2-12\pmod{p}$. This proves the conclusion (ii).

We are done.
\end{proof}

\begin{lemma}\label{lemma34}Suppose that $m=\delta p$, where $\delta \ge 39$, $p\ge 5$ is a prime, $p\not=7$ and $p\nmid \delta$. Then there exist  $1\le a<b\le n$ such that $b^3+b\equiv a^3+a\pmod{m}$.\end{lemma}
\begin{proof}By Lemma \ref{lemma32} and Lemma \ref{lemma33} (ii), we only need to verify $p+\delta L_p\le n$. By \eqref{inequality1}, $n> \frac{\delta p}{3}$ and it suffices to prove  $p+\delta L_p\le \frac{\delta p}{3}$.
Indeed we can prove $\frac{p}{39}+L_p\le \frac{p}{3}$ for all $p\not=7$.
This completes the proof.
\end{proof}

Similarly, we have the following.

\begin{lemma}\label{lemma35}Suppose that $m=\delta p$, where $\delta \ge 13$, $p\ge 165$ is a prime and $p\nmid \delta$. Then there exist  $1\le a<b\le n$ such that $b^3+b\equiv a^3+a\pmod{m}$.\end{lemma}
\begin{proof}It suffices to prove  $p+\delta L_p\le \frac{\delta p}{3}$.
Indeed we can prove $\frac{p}{13}+L_p\le \frac{p}{3}$ for all $p\ge 165$.
This completes the proof.
\end{proof}

\begin{lemma}\label{lemma36}If $p\ge 4000$ and $p\nmid \delta$, then we have
$$\ell_p(\delta)\le \frac{p}{6}.$$\end{lemma}
\begin{proof}We write
\begin{align*}Y=\lfloor\frac{p-1}{6}\rfloor.\end{align*}
It suffices to prove
\begin{align}\label{sumY}\Big|\sum_{1\le x\le Y}\Big(\frac{\delta^2x^2+4}{p}\Big)\Big|<Y-1,\end{align}
since \eqref{sumY} implies that $\Big(\frac{\delta^2x^2+4}{p}\Big)$ can take both $1$ and $-1$ in the range $1\le x\le Y$.
We define
\begin{align*}A=\sum_{-Y\le x\le Y}\Big(\frac{\delta^2x^2+4}{p}\Big).\end{align*}
Note that \eqref{sumY} is equivalent to $|A-1|<2Y-2$, which follows from $|A|<2Y-3$. For $c,x\in \Z$, we have
\begin{align*}\frac{1}{p}\sum_{u=1}^{p}e\big(\frac{u(c-x)}{p}\big)=\begin{cases}
1, \ & \textrm{ if } c\equiv x\pmod{p},
\\ 0, \ & \textrm{ if } c\not\equiv x\pmod{p},\end{cases}\end{align*}
and therefore,
\begin{align*}A=\frac{1}{p}\sum_{1\le u\le p}\sum_{1\le c\le p}\Big(\frac{\delta^2c^2+4}{p}\Big)e(\frac{uc}{p})\sum_{-Y\le x\le Y}e(-\frac{ux}{p}).\end{align*}
By Lemma \ref{lemma31}, we obtain
\begin{align*}|A|\le \frac{2\sqrt{p}}{p}\sum_{1\le u\le p}\Big|\sum_{-Y\le x\le Y}e(-\frac{ux}{p})\Big|,\end{align*}
and by Lemma 4.8 in \cite{YZ} we further have
\begin{align*}|A|\le   2\sqrt{p}(2+\ln p).\end{align*}
The inequality $|A|<2Y-3$ follows from
\begin{align*} 2\sqrt{p}(2+\ln p)<2Y-3.\end{align*}
Note that
\begin{align*}2Y-3>\frac{p}{3}-5.\end{align*}
Now we need to prove
\begin{align}\label{check}2\sqrt{p}(2+\ln p)<\frac{p}{3}-5.\end{align}
It is easy to prove that \eqref{check} holds for $p\ge 4000$.
This completes the proof.
\end{proof}

\begin{lemma}\label{lemma37}Suppose that $m=\delta p$, where $\delta \ge 6$, $p\ge 4000$ is a prime and $p\nmid \delta$. Then there exist $1\le a<b\le n$ such that $b^3+b\equiv a^3+a\pmod{m}$.\end{lemma}
\begin{proof}By Lemma \ref{lemma33} (ii) and Lemma \ref{lemma36}, it suffices to prove  $p+\frac{\delta p}{6}\le \frac{\delta p}{3}$, which holds for $\delta\ge 6$.
This completes the proof.
\end{proof}

\begin{lemma}[Case (i)]\label{lemma38}Let $n\ge 48000$. Suppose that $m=\delta p$, where $\delta \ge 6$, $p\ge 5$ is a prime, $p\not=7$ and $p\nmid \delta$. Then there exist $1\le a<b\le n$ such that $b^3+b\equiv a^3+a\pmod{m}$.\end{lemma}
\begin{proof}In view of Lemma \ref{lemma34}, we only need to consider $\delta<39$. By \eqref{inequality1}, $m\ge 48000$. We deduce that $p=\frac{m}{\delta}\ge \frac{n}{\delta}>165$.
By Lemma \ref{lemma35}, we only need to consider $\delta\le 12$. Now we further have $p=\frac{m}{\delta}\ge \frac{n}{\delta}\ge 4000$ and the desired conclusion follows from Lemma \ref{lemma37}. This completes the proof.
\end{proof}
\begin{remark}\label{remark}One can verify Theorem \ref{theorem1} for $n\le 48000$ with the help of a computer. In fact, Z.-W. Sun has verified the truth of Theorem \ref{theorem1} for $n\le 10^5$. Therefore, the condition $n\ge 48000$ in Lemma \ref{lemma38} can be removed.\end{remark}

\bigskip

\section{The Cases (ii)-(vii)}

\setcounter{lemma}{0}
\setcounter{theorem}{0}
\setcounter{corollary}{0}
\setcounter{equation}{0}

\medskip

The purpose of this section is to deal with cases (ii)-(vii).
\begin{lemma}[Case (ii)]\label{lemma41}Suppose that $m=\delta p^{r}$, where $\delta \ge 4$, $p\ge 5$ is a prime, $r\ge 2$ is a positive integer and $p\nmid \delta$. Then there exist  $1\le a<b\le n$ such that $b^3+b\equiv a^3+a\pmod{m}$.\end{lemma}
\begin{proof} By Lemma \ref{lemma33}, it is sufficient to prove $p^{r}+\delta\frac{p-1}{2}\le n$. By \eqref{inequality1}, $n> \frac{\delta p^{r}}{3}$ and it suffices to prove  $p^{r}+\frac{1}{2}\delta p\le \frac{\delta p^{r}}{3}$. This follows from
\begin{align}\label{ineqin41}(p^{r-1}-\frac{3}{2})(\delta-3)\ge \frac{9}{2}.\end{align}
Since $p^{r-1}-\frac{3}{2}\ge p-\frac{3}{2}\ge \frac{7}{2}$, \eqref{ineqin41} holds if $\delta\ge 5$. In the case $\delta=4$, \eqref{ineqin41} holds if $p^{r-1}\ge 6$. We now only need to consider $\delta=4$, $p=5$, $r=2$, and it is easy to verify that $p^{r}+\delta\frac{p-1}{2}\le \frac{\delta p^{r}}{3}\le n$ holds.

This completes the proof.
\end{proof}

\begin{lemma}[Case (iii)]\label{lemma42}Suppose that $m=2^r$, where $r\in \Z^{+}$. Then there exist $1\le a<b\le n$ such that $b^3+b\equiv a^3+a\pmod{m}$.\end{lemma}
\begin{proof}Note that $2^3+2-1^3-1=2^3$ and $5^3+5-1^3-1=2^7$. For $r\le 3$, we can choose $a=1$ and $b=2$. For $4\le r\le 7$, we can choose $a=1$ and $b=5$.

Now we assume that $r\ge 8$. The proof is the same as that of Lemma 3.7 in \cite{YZ}, and thus we explain it briefly. Since $(a+4)^3+(a+4)-a^3-a=4(3(a+2)^2+5)$, it suffices to find $1\le a\le 2^{r-2}-3$ such that $3(a+2)^2+5\equiv 0\pmod{2^{r-2}}$. For $r\ge 8$, we can find  $3\le x\le 2^{r-2}-1$ such that $3x^2+5\equiv 0\pmod{2^{r-2}}$.  On choosing $a=x-2$, we obtain $3(a+2)^2+5\equiv 0\pmod{2^{r-2}}$. Note that $b=4+a=x+2\le 2^{r-2}+1\le n$. We are done.
\end{proof}

\begin{lemma}[Case (iv)]\label{lemma43}Suppose that $m=2^{r}t$, where $r\ge 2$ is an integer, $t\ge 5$ is an odd number. Then there exist  $1\le a<b\le n$ such that $b^3+b\equiv a^3+a\pmod{m}$.\end{lemma}
\begin{proof}We consider $1\le a\le 2^{r}$ and $b=a+t$. Note that $2^{r}+t\le \frac{2^rt}{3}$ is equivalent to $(2^r-3)(t-3)\ge 9$, which holds  expect that $r=2$ and $t\le 11$. In view of Remark \ref{remark}, we may assume that $n>44$, and by \eqref{inequality1} we have $b\le 2^{r}+t\le n$. It suffices to find $1\le a\le 2^{r}$ such that
$$a^2+a(a+t)+(a+t)^2+1\equiv 0\pmod{2^{r}},$$
and the proof of (3.1) in \cite{YZ} also implies the above conclusion.  We are done.
\end{proof}

\begin{lemma}[Case (v)]\label{lemma44}Suppose that $m=2^r3^s$, where $r,s\in \Z^{+}$. Then there exist  $1\le a<b\le n$ such that $b^3+b\equiv a^3+a\pmod{m}$.\end{lemma}
\begin{proof}By Lemma \ref{lemma43}, we only need to consider either $r=1$ or $s=1$.

We first consider $s=1$. Note that $a^2+a(a+3)+(a+3)^2+1$ is equal to $40=2^3\cdot 5$ when $a=2$. If $r\le 3$, then the desired conclusion follows by choosing $a=2$ and $b=5$. Next we assume $r\ge 4$. Similarly to the proof of (3.1) in \cite{YZ}, we can obtain that for any $j\ge 3$, there exists $1\le a\le 2^{j}-6$ such that
\begin{align*}a^2+a(a+3)+(a+3)^2+1\equiv 0\pmod{2^{j}}.\end{align*}
In particular, there exists $1\le a\le 2^{r}-6$ such that
$a^2+a(a+3)+(a+3)^2+1\equiv 0\pmod{2^{r}}$. The desired conclusion follows by choosing $b=a+3$ and noting that $b\le 2^{r}-3<n$.

Now we consider $r=1$. By \eqref{inequality1}, $n>3^s$. We can choose $a=1$ and $b=1+3^s$.

The proof is complete.
\end{proof}

\begin{lemma}[Case (vi)]\label{lemma45}Suppose that $m=3^{r}\cdot 14$, where $r\in \N$. Then there exist  $1\le a<b\le n$ such that $b^3+b\equiv a^3+a\pmod{m}$.\end{lemma}
\begin{proof}When $m=14$, it suffices to choose $a=1$ and $b=3$. Now we assume that $r\ge 1$. By Lemma \ref{lemma32} (noting that $L_7=3$) and Lemma \ref{lemma33} (ii) (with $p=7$ and $\delta=2\cdot 3^{r}$),
we only need to verify $7+3^{r+1}\cdot 2\le n$. By \eqref{inequality1}, $r=k-3$ and $n>3^{k-1}=3^{r+2}$. Note that $7+3^{r+1}\cdot 2<3^{r+2}$ if $r\ge 1$.
This completes the proof.
\end{proof}

The last task in this section is to consider Case (vii). The proof is as same as that in Section 4 \cite{YZ}. We introduce
\begin{align}X:=X_p=\lfloor\frac{p}{3}\rfloor p^{t-1}.\end{align}

We aim to find $1\le a\not=b\le \frac{n}{\delta}$ such that $\delta^2(a^2+ab+b^2)+1\equiv 0\pmod{p^{t}}$. Then on choosing $a'=\delta a$, $b'=\delta b$, we obtain $a'^3+a'\equiv b'^3+b'\pmod{m}$. By \eqref{inequality1}, we have $X<\frac{n}{\delta}$.

Let
\begin{align}\label{definef}f(a,b)=\delta^2(a^2+ab+b^2).\end{align}Now we introduce
\begin{align}\label{defineN}\mathcal{N}=\sum_{\substack{1\le a,b\le X \\ f(a,b)+1\equiv 0\pmod{p^{t}}}}1\end{align}
and
\begin{align*}\mathcal{N}^{\not=}=\sum_{\substack{1\le a\not=b\le X \\ f(a,b)+1\equiv 0\pmod{p^{t}}}}1.\end{align*}
Note that $\mathcal{N}^{\not=}\ge \mathcal{N}-2$. The main objective is to prove $\mathcal{N}>2$ (and thus $\mathcal{N}^{\not=}>0$.

For $j\ge 1$, we define
\begin{align}\label{defineTj}T_j=\sum_{\substack{1\le c\le p^j \\ (c,p)=1}}\sum_{\substack{1\le a,b\le X }}e\big(\frac{cf(a,b)+c}{p^j}\big).\end{align}

\begin{lemma}\label{lemma46}Let $\mathcal{N}$ and $T_j$ be given in \eqref{defineN} and \eqref{defineTj} respectively. We have
\begin{align*}\mathcal{N}=\frac{X^2}{p^{t}}+\frac{1}{p^{t}}\sum_{j=1}^{t}T_j.\end{align*}
If $1\le j\le t-1$, then
\begin{align*}T_j=X^2p^{-j}\Big(\frac{-3}{p^j}\Big)\mu(p^j),\end{align*}
where $\mu(\cdot)$ is the M\"obius function.
Moreover, we also have
\begin{align*}|T_t|\le 2p^{\frac{3t}{2}}(2+\ln p^t)^2.\end{align*}
\end{lemma}
\begin{proof}The three conclusions are corresponding to Lemma 4.2, Lemma 4.4 and Lemma 4.9 in \cite{YZ} respectively. Although we only considered the case $t=2r$ in \cite{YZ}, both the proofs and the conclusions of Lemmas 4.2, 4.4 and 4.7 in \cite{YZ} are valid for all $t\in \Z^{+}$.\end{proof}

\begin{lemma}\label{lemma47}If $t\ge 2$, then
\begin{align}\label{finalineq1}\mathcal{N}\ge \frac{X^2}{p^{t}}-\Big(\frac{-3}{p}\Big)\frac{X^2}{p^{t+1}}-2p^{t/2}(2+\ln p^{t})^2. \end{align}
If $t=1$, then
\begin{align}\label{finalineq2}\mathcal{N}\ge \frac{X^2}{p^{t}}-2p^{t/2}(2+\ln p^{t})^2. \end{align}
\end{lemma}
\begin{proof}The desired conclusions follow from Lemma \ref{lemma46}.\end{proof}

\begin{lemma}\label{lemma48}Suppose that $m=\delta p^{t}$, where $1\le \delta \le 3$, $p\ge 5$ is a prime, $t\in \Z^{+}$. Suppose further that $p^t\ge 20000^2$. Then we have
\begin{align*}\mathcal{N}^{\not=}>0.\end{align*}
\end{lemma}
\begin{proof}

For $p\ge 7$, we have $\lfloor\frac{p}{3}\rfloor\ge \frac{3}{11}p$ (the equality holds with $p=11$) and $1-\frac{1}{p}\ge \frac{6}{7}$.  Thus for $p\ge 7$, we have
\begin{align}\label{check1}\lfloor\frac{p}{3}\rfloor^2 p^{-2}\Big(1-\big(\frac{-3}{p}\big)\frac{1}{p}\Big)\ge \lfloor\frac{p}{3}\rfloor^2 p^{-2}(1-\frac{1}{p})\ge \frac{3^2\cdot 6}{11^2 \cdot 7}>\frac{6}{125}.\end{align}
For $p=5$, we have
\begin{align}\label{check2}\lfloor\frac{p}{3}\rfloor^2 p^{-2}\Big(1-\big(\frac{-3}{p}\big)\frac{1}{p}\Big)=\frac{6}{125}.\end{align}
We deduce from \eqref{finalineq1}, \eqref{finalineq2}, \eqref{check1} and \eqref{check2} that (for all $p\ge 5$)
\begin{align*}\mathcal{N}\ge \frac{6}{125}p^{t}-2p^{t/2}(2+\ln p^{t})^2. \end{align*}

Since $\mathcal{N}^{\not=}\ge \mathcal{N}-2$, we need to prove
$$\frac{6}{125}p^{t}>2p^{t/2}(2+\ln p^t)^2+2,$$
which follows from
\begin{align}\label{check4}\sqrt{p^t}>\frac{125}{3}(2+\ln p^t)^2+30.\end{align}
On writing $q=\sqrt{p^t}$, our task is to prove $q>\frac{500}{3}(1+\ln q)^2+30$.
Let $g(x)=\sqrt{x-30}-\sqrt{\frac{500}{3}}(1+\ln x)$. Then
$g'(x)=\frac{1}{2\sqrt{x-30}}-\sqrt{\frac{500}{3}}\cdot\frac{1}{x}>\frac{1}{2\sqrt{x}}-\sqrt{\frac{500}{3}}\cdot\frac{1}{x}$ for $x>30$ and $g$ is increasing when $x> \frac{2000}{3}$. Note that $g(20000)>0$. Therefore, $q>\frac{500}{3}(1+\ln q)^2+30$ holds for $q\ge 20000$ and \eqref{check4} holds due to $p^t\ge 20000^2$.
The proof is complete.
\end{proof}

In view of \eqref{inequality1}, for $m=\delta p^t$ (with $1\le \delta \le 3$ and $p\ge 5$ a prime) we have
$$3^{k-1}<n\le \delta p^t<3^{k},$$
and we define
\begin{align}\label{Nstar}\mathcal{N}^\ast=\sum_{\substack{1\le a<b\le 1+3^{k-1} \\ a^3+a\equiv b^3+b\pmod{\delta p^{t}}}}1.\end{align}
We verify $\mathcal{N}^\ast>0$ for $p^t<20000^2$ with the help of a computer.

\begin{lemma}\label{lemma49}Let $\mathcal{N}^\ast$ be given in \eqref{Nstar}. Suppose that $m=\delta p^{t}$, where $1\le \delta \le 3$, $p\ge 5$ is a prime, $t\in \Z^{+}$. Suppose further that $p^t<20000^2$. Then we have
\begin{align*}\mathcal{N}^\ast>0.\end{align*}
\end{lemma}
\begin{proof}This is checked by C++.\end{proof}

\begin{lemma}[Case (vii)]\label{lemma410}Suppose that $m=\delta p^{t}$, where $1\le \delta \le 3$, $p\ge 5$ is a prime, $t\in \Z^{+}$. Then there exist $1\le a<b\le n$ such that $b^3+b\equiv a^3+a\pmod{m}$.\end{lemma}
\begin{proof}The desired conclusion follows from Lemma \ref{lemma48} and Lemma \ref{lemma49}. \end{proof}

\bigskip

\section{The Case (viii)}

\setcounter{lemma}{0}
\setcounter{theorem}{0}
\setcounter{corollary}{0}
\setcounter{equation}{0}

\medskip

It is in Case (viii) that we need to distinguish $n\in \mathcal{E}$ or not in the proof. For $m=3^r\cdot 7$, by \eqref{inequality1} we have $r=k-2$ and
$$n>3^{r+1}.$$

\begin{lemma}\label{lemma51}(i) If $r\equiv 0\pmod{3}$ or $r\equiv 2\pmod{3}$, then we have $\ell_7(3^r)\le 2$.

(ii) If $r\equiv 1\pmod{3}$, then we have $\ell_7(3^r)=3$.
\end{lemma}
\begin{proof}Since $\big(\frac{-3}{7}\big)=1$, $\ell_7(\delta)$ is the smallest positive integer $x$ such that
$\big(\frac{\delta^2x^2+4}{7}\big)=1$.

If $r\equiv 0\pmod{3}$, then $3^{2r}x^2+4\equiv x^2+4\equiv 1\pmod{7}$ for $x=2$ and thus $\ell_7(3^r)\le 2$ (indeed $\ell_7(3^r)=2$ in this case).

If $r\equiv 2\pmod{3}$, then $3^{2r}x^2+4\equiv 4x^2+4\equiv 1\pmod{7}$ for $x=1$ and thus $\ell_7(3^r)=1$.

If $r\equiv 1\pmod{3}$, then $3^{2r}x^2+4\equiv 2x^2+4\pmod{7}$. Note that $\big(\frac{2\cdot 1^2+4}{7}\big)=\big(\frac{2\cdot 2^2+4}{7}\big)=-1$ and $\big(\frac{2\cdot 3^2+4}{7}\big)=1$. Therefore, $\ell_7(3^r)=3$.

This completes the proof.\end{proof}

\begin{lemma}\label{lemma52}Suppose that $m=3^{r}\cdot 7$, where $r\equiv 0\pmod{3}$ or $r\equiv 2\pmod{3}$. Then there exist  $1\le a<b\le n$ such that $b^3+b\equiv a^3+a\pmod{m}$.\end{lemma}
\begin{proof}For $r=0$, we can choose $a=1$ and $b=3$. By Lemma \ref{lemma51} (i), $\ell_7(3^r)\le 2$. For $r\ge 2$, the desired the conclusion follows from Lemma \ref{lemma33} (ii) on noting that $7+3^r\cdot 2\le 3^{r+1}<n$.\end{proof}

\begin{lemma}\label{lemma53}Suppose that $m=3^{r}\cdot 7$, where $r\equiv 1\pmod{3}$. Suppose further that either $r\equiv 1\pmod{6}$ or $n\not\in \mathcal{E}$. Then there exist  $1\le a<b\le n$ such that $b^3+b\equiv a^3+a\pmod{m}$.\end{lemma}
\begin{proof}Note that $a^2+a(a+3^{r+1})+(a+3^{r+1})^2+1=3a^2+3^{r+2}a+3^{2r+2}+1$. It suffices to find $a\in \Z^{+}$ such that $3a^2+3^{r+2}a+3^{2r+2}+1\equiv 0\pmod{7}$ and $a+3^{r+1}\le n$. Note that for $r\equiv 1\pmod{3}$, we have $3^{2r+2}\equiv 4\pmod{7}$. On writing $r=6s+1+3t$ with $s\in \N$ and $t\in \{0,1\}$, we have
\begin{align}\label{conginlemma53}3a^2+3^{r+2}a+3^{2r+2}+1 \equiv 3a^2+3^{3t+3}a+5\equiv 3a^2+(-1)^{t+1}a+5\pmod{7}.\end{align}

If $t=0$, then by \eqref{conginlemma53} we can choose $a=1$ such that $3a^2+3^{r+2}a+3^{2r+2}+1\equiv 0\pmod{7}$ and $a+3^{r+1}=1+3^{r+1}\le n$.

If $t=1$, then $n\not\in\mathcal{E}$ and $n\ge 3^{r+1}+3$. By \eqref{conginlemma53}, we choose $a=3$ such that $3a^2+3^{r+2}a+3^{2r+2}+1\equiv 0\pmod{7}$. Note that
$b=3+3^{r+1}\le n$. We are done.
\end{proof}

\begin{lemma}[Case (viii), Part 1]\label{lemma54}Let $n\not\in \mathcal{E}$. Suppose that $m=3^{r}\cdot 7$. Then there exist  $1\le a<b\le n$ such that $b^3+b\equiv a^3+a\pmod{m}$.\end{lemma}
\begin{proof}The desired conclusion follows from Lemmas \ref{lemma52}-\ref{lemma53}.\end{proof}

\begin{lemma}[Case (viii), Part 2]\label{lemma55}Let $n=3^{6s+5}+1$ or $n=3^{6s+5}+2$. Suppose that $m=3^{r}\cdot 7$. Then $a^3+a(1\le a\le n)$ are pairwise distinct modulo $m$.\end{lemma}
\begin{proof}By \eqref{inequality1}, we have $r=6s+4$. Suppose otherwise that we can find $1\le a<b\le n$ such that $b^3+b\equiv a^3+a\pmod{m}$. Note that $3\nmid (a^2+ab+b^2+1)$ for any $a,b\in \Z$, and we conclude that  $b=a+3^{6s+4}c$ for some $c\in \Z^{+}$. Since $b=a+3^{3s+4}c<n$, we have $c\le 3$. Write $\delta=3^{6s+4}$. Then $b^3+b\equiv a^3+a\pmod{m}$ implies that $a^2+ab+b^2+1\equiv 0\pmod{7}$, which is equivalent to
\begin{align*}(6a+3\delta c)^2\equiv -3\delta^2c^2-12\pmod{7}.\end{align*}
Therefore, we have $\Big(\frac{-3\delta^2c^2-12}{7}\Big)\in\{0,1\}$. By Lemma \ref{lemma51} (ii), $\ell_7(\delta)=3$ for $\delta=3^{6s+4}$, and we obtain $c\ge \ell_7(\delta)=3$. Now we conclude that $c=3$ and $b=a+3^{6s+5}$.
 Then we deduce that $a^2+ab+b^2+1\equiv 3a^2+a+5\equiv 0\pmod{7}$, and which implies $a\ge 3$ and $b=a+3^{3s+5}\ge 3+3^{3s+5}$. This is a contradiction to $b\le n$. The proof is complete.
\end{proof}

\bigskip

\noindent {\it Proof of Lemmas \ref{lemma21}-\ref{lemma22}.} In view of Lemma \ref{lemma38}, Remark \ref{remark}, Lemma \ref{lemma41}, Lemma \ref{lemma42}, Lemma \ref{lemma43}, Lemma \ref{lemma44}, Lemma \ref{lemma45}, Lemma \ref{lemma410}, Lemma \ref{lemma54} and Lemma \ref{lemma55}, we only need to prove that each positive integer $m$ restricted by \eqref{inequality1} must satisfy (at least) one of the $8$ cases in Section 2. By \eqref{inequality1}, $m$ is not a power of $3$.

If $m$ has no prime factors greater than $3$, then $m$ belongs to Case (iii) or (v). Next we assume that $m$ has two distinct prime factors greater than $3$. We write $m=m'p_1^{r_1}p_2^{r_2}$, where $p_1\not =p_2$ are two primes, $p_1\nmid m'$, $p_2\nmid m'$ and $r_1,r_2\in \Z^{+}$. Without loss of generality, we further assume that $p_1\not=7$ and $p_2\ge 7$. Let $\delta=m'p_2^{r_2}$. Then $m=\delta p_1^{r_1}$, $\delta\ge 7$ and $p_1\nmid \delta$. We can see that $m$ belongs to either Case (i) or Case (ii).

Now we assume that $m$ has only one prime factor greater than $3$. We write $m=2^{i}3^{j}p^{r}$, where $p\ge 5$ is a prime, $r\in \Z^{+}$ and $i,j\in \N$. Note that if $i\ge 2$, then $m$ satisfies the condition of Case (iv). We discuss $i=1$ and $i=0$ below.

We first consider $i=1$. If $j=0$, then $m$ belongs to Case (vii). If $j\ge 1$ and $r\ge 2$, then $m$ satisfies the condition of Case (ii). If $j\ge 1$, $r=1$ and $p\not=7$, then $m$ satisfies the condition of Case (i). If $j\ge 1$, $r=1$ and $p=7$, then $m$ satisfies the condition of Case (vi).

Now we consider $i=0$. If $0\le j\le 1$, then $m$ belongs to Case (vii). If $j\ge 2$ and $r\ge 2$, then $m$ satisfies the condition of Case (ii). If $j\ge 2$, $r=1$ and $p\not=7$, then $m$ satisfies the condition of Case (i). If $j\ge 2$, $r=1$ and $p=7$, then $m$ satisfies the condition of Case (viii).

We have proved that $m$ subject to \eqref{inequality1} must satisfy (at least) one of the $8$ cases in Section 2.

According to the remark before Lemma \ref{lemma21}, we also complete the proof of Theorem \ref{theorem1}.

\end{document}